\documentclass[12pt]{article}
\usepackage{amsmath,amsfonts,amssymb,amsthm}
\usepackage{graphics}
\usepackage{epsfig}
\usepackage{esint}
\usepackage{color}
\newcommand{\ignore}[1]{}

\newcommand{\e}{\varepsilon}
\newcommand{\R}{\mathbb{R}}
\newcommand{\ho}{{hom}}

\newcommand{\expec}[1]{\langle #1 \rangle}
\newcommand{\Expec}[1]{\left\langle #1 \right\rangle}
\newtheorem{definition}{Definition}
\newtheorem{proposition}{Proposition}
\newtheorem{theorem}{Theorem}
\newtheorem{remark}{Remark}
\newtheorem{lemma}{Lemma}

\title{Quantitative estimates on the periodic approximation of the corrector in stochastic homogenization}
\author{A. Gloria
\thanks{Universit\'e Libre de Bruxelles (ULB), Brussels, Belgium \&
  Team MEPHYSTO, Inria Lille - Nord Europe, Villeneuve d'Ascq,
  France, \textrm{agloria@ulb.ac.be}}
\and F. Otto\thanks{Max Planck Institute for Mathematics in the Sciences, Inselstr. 22, 04103 Leipzig,
Germany, \textrm{felix.otto@mis.mpg.de}}}
\begin{document}
\maketitle


\section{Introduction, setting, and main results}

\subsection{Introduction}

This article corresponds to a course given by F. Otto at the summer school CEMRACS 2013 in Luminy, France.

\medskip

The qualitative theory of stochastic homogenization dates back to the seminal contributions of Papanicolaou and Varadhan \cite{Papanicolaou-Varadhan-79}, and Kozlov \cite{Kozlov-79}. Let $D$ a bounded domain of $\R^d$, $f\in H^{-1}(D)$, and $A$ be a stationary ergodic random field. Then the weak solution $u_\e\in H^1_0(D)$ of
$$
-\nabla \cdot a(\frac{\cdot}{\e})\nabla u_\e\,=\,f
$$
almost surely weakly converges in $H^1(D)$ to the unique weak solution $u_0\in H^1_0(D)$ of
\begin{equation}\label{eq:hom}
-\nabla \cdot a_\ho\nabla u_0\,=\,f,
\end{equation}
where $a_\ho$ is a deterministic matrix (independent of $f$) called the homogenized matrix.
As far as quantitative estimates are concerned there are only few contributions in the literature.
A first general comment is that ergodicity alone is not enough to obtain convergence rates, so that 
mixing properties have to be assumed on the random coefficients $a$.
Besides the optimal estimates in the one-dimensional case by Bourgeat and Piatnitskii \cite{Bourgeat-99},
the first and still unsurpassed contribution in the linear case is due to Yurinski{\u\i} who proved in \cite[(0.10)]{Yurinskii-86} that for $d>2$ and for mixing coefficients with an algebraic decay (not necessarily integrable), there exists $\gamma>0$
such that 
\begin{equation}\label{eq:intro-3}
\expec{|u_\e-u_\ho|^2}\,\lesssim \,\e^{\gamma}.
\end{equation}

\medskip

In the case of discrete elliptic equations, the quantitative theory is much more developed and started with the inspiring unpublished work by Naddaf and Spencer \cite{Naddaf-Spencer-98}, who understood the right notion of mixing condition for stochastic homogenization in the form of a spectral gap estimate for the Glauber dynamics. In \cite{Gloria-Otto-09,Gloria-Otto-09b} we proved the first optimal quantitative estimates on the corrector, on the variance of the spatial averages of the energy density of the corrector, and estimates on the difference between the homogenized coefficients and their approximations using a massive term in the corrector equation. In collaboration with Neukamm in \cite{Gloria-Neukamm-Otto-14} we then obtained optimal estimates in any dimension on the corrector, on its periodic and massive approximations, and on the approximations of the homogenized coefficients by periodization or using a massive term. Marahrens and the second author proved in \cite{Marahrens-Otto-13} optimal annealed estimates on the Green function, which in turn allowed us and Neukamm to prove an optimal quantitative two-scale expansion of $u_\e$ for all $d\ge 2$, thus improving considerably  
\eqref{eq:intro-3} in the discrete setting.
More recent results go beyond variance estimates and address statistical properties of solutions such as central limit theorems (CLT). In particular, for the approximation of homogenized coefficients, a CLT was recently obtained by Biskup, Salvi and Wolff \cite{Biskup-Salvi-Wolff-14} in the case of small ellipticity ratio (that is, for a coefficient field $a$ close to identity)
 for independent and identically distributed coefficients, see also Rossignol \cite{Rossignol-12}. An optimal quantitative version of this CLT (which estimates in particular the Wasserstein distance to normality) is proved by Nolen and the first author in \cite{Gloria-Nolen-14} (without the smallness assumption on the ellipticity ratio). The route towards a CLT for the solution $u_\e$ itself started for $d>2$ with the contribution by Mourrat and the second author on the structure of the asymptotic covariance of the corrector, see \cite{Mourrat-Otto-14}.

\medskip

The present contribution concerns the case of continuum linear elliptic equations.
It is based on \cite{Gloria-Otto-10b}, which is a complete continuum version of \cite{Gloria-Otto-09,Gloria-Otto-09b} (with in addition optimal results for $d=2$). In the present contribution we establish quantitative results on the periodic approximation of the corrector equation for the stochastic homogenization of linear elliptic equations in divergence form, 
when the diffusion coefficients satisfy a spectral gap estimate in probability, and for $d>2$.
The main difference with respect to the first part of \cite{Gloria-Otto-10b} is that we avoid here the use of Green's functions and more directly rely on the De Giorgi-Nash-Moser theory. Let us also mention the work \cite{Nolen-11} by Nolen, which follows \cite{Gloria-Otto-09} to establish estimates on the corrector, and proves a normal approximation result for a class of random coefficients in this continuum setting.

\subsection{Setting}

We start by introducing the relevant deterministic notions: 
The corrector $\phi(a;\cdot)$
and the homogenized coefficient 
$a_{hom}^L(a)$ 
for an arbitrary coefficient field $a$ on the torus of side length $L$,
which we sometimes denote as $[-\frac{L}{2},\frac{L}{2})^d$ to single
out the point $0$. 

\begin{definition}\label{defi:hom-coeff}
\mbox{}\\
{\sc Space of coefficient fields}. For a given side-length $L$ let
$\Omega$ be the space of all $[-\frac{L}{2},\frac{L}{2})^d$-periodic 
fields of $d\times d$ matrices $a$  
that are
uniformly elliptic in the sense
\begin{equation}\nonumber
\forall\;x\in\Big[-\frac{L}{2},\frac{L}{2}\Big)^d,\;\xi\in\mathbb{R}^d\quad\lambda|\xi|^2\le \xi\cdot a(x)\xi, \ |a(x)\xi|\le|\xi|,
\end{equation}
where $\lambda>0$ is a number fixed throughout the article.

\smallskip

{\sc Corrector}. For given $a\in\Omega$, the corrector $\phi(a;\cdot)$ is an 
$[-\frac{L}{2},\frac{L}{2})^d$-periodic 
function defined through the elliptic equation
\begin{equation}\label{T.3}
-\nabla\cdot a(\nabla\phi(a;\cdot)+\xi)=0\quad\mbox{and}\quad\int_{[-\frac{L}{2},\frac{L}{2})^d}\phi(a;\cdot)=0,
\end{equation}
where $\xi\in\mathbb{R}^d$ with $|\xi|=1$ is a direction which is fixed throughout the article.
For further reference we note that $\phi$ is ``stationary'' in the sense of
\begin{equation}\label{L2.11}
\phi(a(\cdot+z),x)=\phi(a,x+z)
\end{equation}
for all points $x\in\mathbb{R}^d$, coefficient fields $a\in\Omega$, and shift vectors $z\in\mathbb{R}^d$.

\smallskip

{\sc Homogenized coefficient}. The homogenized coefficient in directions
$\xi,\xi'$ is defined via
\begin{equation}\label{T.2}
\xi'\cdot a_{hom}^L(a)\xi:=L^{-d}\int_{[-\frac{L}{2},\frac{L}{2})^d}
\xi'\cdot a(\nabla\phi(a;\cdot)+\xi),
\end{equation}
where $\xi'$ with $|\xi'|=1$ is a direction which is fixed throughout the article.
\end{definition}
We now introduce our example of an ensemble on the space of coefficient fields on the torus. 

\smallskip

\begin{definition}
By the ``Poisson ensemble'' we understand the following probability measure on $\Omega$:

\smallskip

Let the configuration of points $X:=\{X_n\}_{n=1,\cdots,N}$ on the torus
be distributed according to the Poisson
point process with density one. This means the following
\begin{itemize}
\item 
For any two disjoint (Lebesgue measurable) subsets $D$ and $D'$ of the torus we have that
the configuration of points in $D$ and the configuration of points in $D'$ are independent.
In other words, if $\zeta$ is a function of $X$ that depends on $X$ only through $X_{|D}$
and $\zeta'$ is a function of $X$ that depends on $X$ only through $X_{|D'}$ we have
\begin{equation}\label{L4.23}
\Expec{\zeta\zeta'}_0=\Expec{\zeta}_0\Expec{\zeta'}_0,
\end{equation}
where $\Expec{\cdot}_0$ denotes the expectation w.\ r.\ t.\ the Poisson point process.
\item For any (Lebesgue measurable) subset $D$ of the torus, 
the number of points in $D$ is Poisson distributed;
the expected number is given by the Lebesgue measure of $D$.
\end{itemize}
Note that $N$ is random, too.

\smallskip

With any realization $X=\{X_n\}_{n=1,\cdots,N}$ of the Poisson point process, 
we associate the coefficient field $a\in\Omega$ via
\begin{equation}\label{X}
a(x)=\left\{\begin{array}{ccc}
\lambda&\mbox{if}&x\in\bigcup_{n=1}^{N}B_1(X_n)\\
1&\mbox{else}\end{array}\right\}{\rm id}.
\end{equation}
Here and throughout the article, balls like $B_1(X_n)$ refer to the distance function
of the torus.
This defines a probability measure on $\Omega$ by ``push-forward'' of $\Expec{\cdot}_0$.
We denote the expectation w.\ r.\ t.\ this ensemble with $\Expec{\cdot}$.
\end{definition}
\begin{remark}
The coefficient field $a$ associated with the Poisson point process is a stationary ergodic coefficient field to which the theory by Papanicolaou and Varadhan applies. In particular, as a direct consequence of the homogenization result, we have almost surely 
$$
\lim_{L\uparrow \infty} a_{hom}^L(a)\,=\,a_{hom},
$$
where $a_{hom}$ is as in \eqref{eq:hom} and $a_{hom}^L(a)$ is as in Definition~\ref{defi:hom-coeff}.
The random variable $a_{hom}^L(a)$ is called the approximation by periodization of the homogenized coefficient $a_{hom}$.
\end{remark}
For our result, we only need the following two properties of the Poisson ensemble.

\begin{lemma}\label{L4}
\mbox{}\\
{\sc Stationarity}.
The Poisson ensemble is
stationary which means that for any shift vector $z\in\mathbb{Z}^d$
the random field $a$ and its shifted version $a(\cdot+z)\colon x\mapsto a(x+z)$ have the same distribution.
In other words, for any (integrable) function $\zeta\colon\Omega\rightarrow\mathbb{R}$ 
(which we think of as a random variable) we have
that $a\mapsto \zeta(a(\cdot+z))$ and $\zeta$ have the same expectation:
\begin{equation}\label{stat}
\Expec{\zeta(a(\cdot+z))}=\Expec{\zeta}.
\end{equation}

\smallskip

{\sc Spectral Gap Estimate}.
The Poisson ensemble satisfies a Spectral Gap Estimate
by which we understand the following: There exists a radius $R$ only depending on $d$
such that for any function $\zeta\colon\Omega\rightarrow\mathbb{R}$, we have
\begin{equation}\label{L3.10}
\Expec{(\zeta-\Expec{\zeta})^2}\le
\Expec{\sum_{z\in\mathbb{Z}^d\cap[-\frac{L}{2},\frac{L}{2})^d}({\rm osc}_{B_R(z)}\zeta)^2}.
\end{equation}
Here, for a (Lebesgue measurable) subset $D$ of the torus, the 
(essential) oscillation ${\rm osc}_{D}\zeta$ of $\zeta$
with respect to $D$ is a random variable defined through
\begin{eqnarray}\nonumber
({\rm osc}_D\zeta)(a)&=&\sup\{\zeta(\tilde a)|\tilde a\in\Omega\;\mbox{with}\;\tilde a=a\;\mbox{outside}\;D\}\nonumber\\
                   & &-\inf\{\zeta(\tilde a)|\tilde a\in\Omega\;\mbox{with}\;\tilde a=a\;\mbox{outside}\;D\}.\label{L4.21}
\end{eqnarray}
It measures how sensitively $\zeta(a)$ depends on $a_{|D}$. Note that $({\rm osc}_D\zeta)(a)$ does
not depend on $a_{|D}$.
\end{lemma}

\subsection{Main results}

The main result of the article is a Central Limit Theorem-type scaling of the
variance of the homogenized coefficient in terms of the system volume $L^d$.

\begin{theorem}\label{T}
Suppose $\Expec{\cdot}$ is stationary and satisfies the Spectral Gap Estimate.
Then we have the following estimate on the variance of the periodic approximation of the homogenized coefficient
\begin{equation}\nonumber
\Expec{(\xi'\cdot a_{hom}^L\xi-\Expec{\xi'\cdot a_{hom}^L\xi})^2}
\le C(d,\lambda) L^{-d}.
\end{equation}
\end{theorem}

In this article, we prove Theorem \ref{T} only for $d>2$. We shall derive it from
the following result of independent interest, which is only true for $d>2$.

\begin{proposition}\label{P}
Let $d>2$ and suppose $\Expec{\cdot}$ is stationary and satisfies the Spectral Gap Estimate.
Then all moments of the corrector are bounded independently of $L$, that is, for any $1\le p<\infty$ we have
\begin{equation}\nonumber
\Expec{\phi^{2p}}\le C(d,\lambda,p).
\end{equation}
Here and in the entire text, we write $\phi^{2p}$ for $(\phi^2)^p$, so that expressions
like above make sense also for a non-integer exponent $p$.

\end{proposition}


\section{Auxiliary results}

We need the following $L^p(\Omega)$-version of the Spectral Gap Estimate.

\begin{lemma}\label{L3}
Let $\Expec{\cdot}$ satisfy the Spectral Gap Estimate. Then it satisfies an
$L^p(\Omega)$-version of a Spectral Gap Estimate in the following sense: 
Let $R$ be the radius from (\ref{L3.10}). Then we have for any ($2p$-integrable) function 
$\zeta\colon\Omega\rightarrow\mathbb{R}$ and any $1\le p<\infty$ 
\begin{equation}\label{L3.b}
\Expec{(\zeta-\Expec{\zeta})^{2p}}\lesssim
\Expec{\Big(\sum_{z\in\mathbb{Z}^d\cap[-\frac{L}{2},\frac{L}{2})^d}({\rm osc}_{B_R(z)}\zeta)^2\Big)^p}.
\end{equation}
Here $\lesssim$ means up to a generic constant that only depends on $p$.

\end{lemma}

Together with the previous lemma, the following lemma gives an estimate of
$\phi$ in terms of $\nabla\phi+\xi$.

\begin{lemma}\label{L2}
Suppose $d>2$ and that $\Expec{\cdot}$ is stationary. Then we have for any 
$\frac{d}{d-2}< p<\infty$ and any $R$ large enough
\begin{equation}\label{L2.b}
\Expec{\Big(\sum_{z\in\mathbb{Z}^d\cap[-\frac{L}{2},\frac{L}{2})^d}({\rm osc}_{B_R(z)}\phi(\cdot;0))^2\Big)^p}
\lesssim\Expec{\Big(\int_{B_1}|\nabla\phi+\xi|^2\Big)^p},
\end{equation}
where $\lesssim$ means up to a generic constant only depending on $d$, $\lambda$, $p$, and $R$.
\end{lemma}

\medskip

The last lemma in turn gives an estimate of $\nabla\phi+\xi$ in terms of $\phi$.

\begin{lemma}\label{L1}
Suppose that $\Expec{\cdot}$ is stationary. Then we have for any $2\le p<\infty$
\begin{equation}\nonumber
\Expec{\Big(\int_{B_1}|\nabla\phi+\xi|^2\Big)^p}\lesssim\Expec{\phi^{2(p-1)}}+1,
\end{equation}
where $\lesssim$ means up to a generic constant only depending on $d$, $\lambda$, and $p$.
\end{lemma}


\section{Proofs}

{\sc Proof of Lemma \ref{L4}}.

{\bf Step 1}. Generalization and reduction. The most natural form of the result of the lemma is
the following: For any measurable partition $D_1,\cdots,D_N$ of the torus we have
\begin{equation}\label{L4.4}
\Expec{(\zeta-\Expec{\zeta})^2}\le\Expec{\sum_{n=1}^N({\rm osc}_{B_1(D_n)}\zeta)^2},
\end{equation}
where $B_1(D)$ is the set of all points on the torus that have distance less than one to $D$. In this
step, we will derive this from the following similar estimate on the Poisson point process itself:
\begin{equation}\label{L4.20}
\Expec{(\zeta_0-\Expec{\zeta_0}_0)^2}_0\le\Expec{\sum_{n=1}^N({\rm osc}_{0,D_n}\zeta_0)^2}_0.
\end{equation}
Here $\Expec{\cdot}_0$ denotes the expectation w.\ r.\ t.\ to the Poisson point process
$X:=\{X_n\}_{n=1,\cdots,N}$,
$\zeta$ is a (square integrable) function of the point configuration $X$,
and the oscillation ${\rm osc}_{0,D}$ is defined in a similar way to (\ref{L4.21}):
\begin{eqnarray}
({\rm osc}_{0,D}\zeta_0)(X)&=&\sup\{\zeta_0(\tilde X)|\;\tilde X=X\mbox{outside}\;D\}\nonumber\\
                   & &-\inf\{\zeta_0(\tilde X)|\;\tilde X=X\mbox{outside}\;D\}.\label{osc0}
\end{eqnarray}
Indeed, (\ref{L4.4}) is an immediate consequence of (\ref{L4.20}) because of
the following two facts.
\begin{itemize}
\item We recall that (\ref{X}) defines a mapping $X\mapsto a$ from point configurations
to coefficient fields. As such, it pulls back functions
according to $\zeta_0(X)=\zeta(a(X))$ and pushes forward the ensemble according to
\begin{equation}\label{L4.22}
\Expec{\zeta}=\Expec{\zeta_0}_0.
\end{equation}
In particular we have for the variance
\begin{equation}\nonumber
\Expec{(\zeta-\Expec{\zeta})^2}
=\Expec{(\zeta_0-\Expec{\zeta_0}_0)^2}_0.
\end{equation}
\item By definition (\ref{X}), if the point configurations $X$ and $\tilde X$ coincide outside  $D$,
then the corresponding coefficient fields $a(X;\cdot)$ and $a(\tilde X;\cdot)$ coincide outside $B_1(D)$.
Hence for a given configuration $X$, the set
$\{a(\tilde X)|\tilde X=X\;\mbox{outside}\;D\}$ is contained
in the set $\{\tilde a|\tilde a=a(X)$ $\;\mbox{outside}\;B_1(D)\}$ so that
\begin{eqnarray*}
\lefteqn{\sup\{\zeta(a(\tilde X))|\tilde X=X\;\mbox{outside}\;D\}}\\
&\le&\sup\{\zeta(\tilde a)|\tilde a=a(X)\;\mbox{outside}\;B_1(D)\},
\end{eqnarray*}
and the opposite inequality if we replace the supremum by the infimum.
From the definitions (\ref{L4.21}) and (\ref{osc0}) of the oscillation we thus see
\begin{equation}\nonumber
({\rm osc}_{0,D}\zeta_0)(X)\le({\rm osc}_{B_1(D)}\zeta)(a(X)).
\end{equation}
By (\ref{L4.22}), this implies as desired
\begin{equation}\nonumber
\Expec{({\rm osc}_{0,D}\zeta_0)^2}_0\le
\Expec{({\rm osc}_{B_1(D)}\zeta)^2}.
\end{equation}
\end{itemize}

\medskip

{\bf Step 2}. Conditional expectations and independence. From now on, we prove statement (\ref{L4.20})
for the Poisson point process. For brevity, we drop the subscript $0$.

\smallskip

For a given (Lebesgue measurable) subset $D$ of the torus,
we denote by $\Expec{\cdot|D}$ the expectation conditioned on the restriction $X_{|D}$ of the
(random) point configuration $X$ on $D$. We note that for a function
$\zeta\colon\Omega\rightarrow\mathbb{R}$ which is square integrable, $\Expec{\zeta|D}$ is the $L^2(\Omega)$-ortho\-gonal
projection of $\zeta$ onto the space of square integrable functions $\tilde\zeta\colon\Omega\rightarrow\mathbb{R}$
that only depend on $X$ via $X_{|D}$. 

\smallskip

With help of these conditional expectations the 
independence assumption (\ref{L4.23}) can be rephrased as follows: For any
two (Lebesgue measurable) subsets $D$, $D'$ that are {\it disjoint}
and any (square integrable) function $\zeta\colon\Omega\rightarrow\mathbb{R}$ 
that does {\it not} depend on $X_{|D}$ we have
\begin{equation}\label{L4.11}
\Expec{\zeta|D\cup D'}=\Expec{\zeta|D'}.
\end{equation}
Here comes the argument: By definition of conditional expectation, (\ref{L4.11}) follows if for
any pair of (bounded and measurable) test functions $u$ and $u'$ which only depend
on $X_{|D}$ and $X_{|D'}$, respectively, we have
\begin{equation}\nonumber
\Expec{\zeta u u'}=\Expec{\Expec{\zeta|D'} u u'}.
\end{equation}
Indeed, on the one hand, since $\zeta$ only depends on $X_{|D^c}$ (where $D^c$ denotes the complement
of $D$) and $u'$ does only depend on $X_{|D'}$ (and thus a fortiori only on $X_{|D^c}$) while $u$
only depends on $X_{|D}$, we have from (\ref{L4.23}):
\begin{equation}\nonumber
\Expec{\zeta u u'}=\Expec{\zeta u'}\Expec{ u}.
\end{equation}
On the other hand, since $\Expec{\zeta|D'} u'$ only depends on $X_{|D'}$ (and in particular
only on $X_{|D^c}$) while $u$ only depends on $X_{|D}$, we have from (\ref{L4.23}):
\begin{equation}\nonumber
\Expec{\Expec{\zeta|D'} u u'}=\Expec{\Expec{\zeta|D'} u'}\Expec{ u}
=\Expec{\Expec{\zeta u'|D'}}\Expec{ u}=\Expec{\zeta u'}\Expec{ u},
\end{equation}
where the middle identity holds since $u'$ only depends on $X_{|D'}$.

\medskip

{\bf Step 3}. Conditional expectation and oscillation. For any (Lebesgue measurable) disjoint subsets $D$ and $D'$ of the torus
and any (square integrable) function, we have
\begin{equation}
|\Expec{\zeta|D\cup D'}-\Expec{\zeta|D'}|\le\Expec{{\rm osc}_{D}\zeta|D'}.\label{L4.8}
\end{equation}
By exchanging $\zeta$ with $-\zeta$, we see that it is enough to show
\begin{equation}\label{L4.25}
\Expec{\zeta|D\cup D'}\le\Expec{\zeta|D'}+\Expec{{\rm osc}_D\zeta|D'},
\end{equation}
We note that $\sup_D\zeta\le\zeta+{\rm osc}_D\zeta$, where we've set for abbreviation
\begin{equation}\nonumber
(\sup_{D}\zeta)(X):=\sup\{\zeta(\tilde X)|\tilde X=X\;\mbox{outside}\;D\}.
\end{equation}
Hence (\ref{L4.25}) follows from
\begin{equation}\nonumber
\Expec{\zeta|D\cup D'}\le\Expec{\sup_D\zeta|D'}.
\end{equation}
The latter inequality can be seen as follows
\begin{eqnarray*}
\lefteqn{\Expec{\zeta|D\cup D'}}\nonumber\\
&\le&\Expec{\sup_D\zeta|D\cup D'}\quad\mbox{since}\;\zeta\le\sup_D\zeta\\
&\stackrel{(\ref{L4.11})}{=}&
\Expec{\sup_D\zeta|D'}\quad\mbox{since}\;\sup_D\zeta\;\mbox{does not depend on }\;X_{|D}.
\end{eqnarray*}

{\bf Step 4}. Martingale decomposition. For conciseness, we only prove (\ref{L4.20})
for $N=3$. So let $\{D_1,D_2,D_3\}$ be a partition of the torus, we claim
\begin{eqnarray}\label{L4.12}
\lefteqn{\Expec{(\zeta-\Expec{\zeta})^2}}\\
&=&\Expec{(\zeta-\Expec{\zeta|D_1\cup D_2})^2}
+\Expec{(\Expec{\zeta|D_1\cup D_2}-\Expec{\zeta|D_1})^2}
+\Expec{(\Expec{\zeta|D_1}-\Expec{\zeta})^2}.\nonumber
\end{eqnarray}
Indeed, this follows from the fact that
\begin{equation}\nonumber
\zeta-\Expec{\zeta|D_1\cup D_2},\,
\Expec{\zeta|D_1\cup D_2}-\Expec{\zeta|D_1},\,
\Expec{\zeta|D_1}-\Expec{\zeta}\mbox{ are } L^2(\Omega)-\mbox{ortho\-gonal}.
\end{equation}
The latter can be seen as follows:
By definition of $\Expec{\cdot|D}$
as $L^2(\Omega)$-orthogonal projection, the two last functions $\Expec{\zeta|D_1\cup D_2}-\Expec{\zeta|D_1}$ and
$\Expec{\zeta|D_1}-\Expec{\zeta}$ do only depend on $X_{|D_1\cup D_2}$, so that they
are orthogonal to the first function $\zeta-\Expec{\zeta|D_1\cup D_2}$. It remains
to argue that the two last functions $\Expec{\zeta|D_1\cup D_2}-\Expec{\zeta|D_1}$ and
$\Expec{\zeta|D_1}-\Expec{\zeta}$ are orthogonal. To that purpose, we rewrite the middle function
as
\begin{equation}\nonumber
\Expec{\zeta|D_1\cup D_2}-\Expec{\zeta|D_1}
=\zeta'-\Expec{\zeta'|D_1}\quad\mbox{where}\;\zeta':=\Expec{\zeta|D_1\cup D_2}.
\end{equation}
Since the last function only depends on $X_{|D_1}$, they are orthogonal.

\medskip

{\bf Step 5}. Conclusion, i.\ e.\ (\ref{L4.4}) for $N=3$. By Step 4, it remains to estimate
the three r.\ h.\ s.\ terms of (\ref{L4.12}). For the first term, we use (\ref{L4.8}) 
with $D'=D_1\cup D_2$ and $D=D_3$
and obtain because of $\zeta=\Expec{\zeta|D_1\cup D_2\cup D_3}$
\begin{eqnarray*}
\Expec{(\zeta-\Expec{\zeta|D_1\cup D_2})^2}&\le&
\Expec{\Expec{{\rm osc}_{D_3}\zeta|D_1\cup D_2}^2}\\
&\stackrel{\mbox{\small Jensen}}{\le}&
\Expec{\Expec{({\rm osc}_{D_3}\zeta)^2|D_1\cup D_2}}
=
\Expec{({\rm osc}_{D_3}\zeta)^2}.
\end{eqnarray*}
The other two terms follow the same way.

\bigskip


{\sc Proof of Lemma \ref{L3}}.

W.\ l.\ o.\ g.\ we may assume that $\Expec{\zeta}=0$.

{\bf Step 1}. Application of the original Spectral Gap Estimate to $\zeta^p$. We claim that this yields
\begin{equation}\label{L3.11}
\Expec{\zeta^{2p}}\lesssim\Expec{\zeta^p}^2
+\Expec{\Big(\sum_{z\in\mathbb{Z}^d\cap[-\frac{L}{2},\frac{L}{2})^d}({\rm osc}_{B_R(z)}\zeta)^2\Big)^p}.
\end{equation}
Indeed, (\ref{L3.10}) applied to $\zeta^p$ at first gives
\begin{equation}\label{L3.13}
\Expec{(\zeta^{p}-\Expec{\zeta^p})^2}\lesssim
\Expec{\sum_{z\in\mathbb{Z}^d\cap[-\frac{L}{2},\frac{L}{2})^d}({\rm osc}_{B_R(z)}(\zeta^p))^2}.
\end{equation}
Using the triangle inequality in $L^2(\Omega)$ on the l.\ h.\ s.\ of (\ref{L3.13}) in form of
\begin{equation}\nonumber
\Expec{(\zeta^p)^2}^\frac{1}{2}\le\Expec{(\zeta^p-\Expec{\zeta^p})^2}^\frac{1}{2}
+|\Expec{\zeta^p}|,
\end{equation}
we see that (\ref{L3.11}) follows from (\ref{L3.13}) by Young's inequality provided we can show
\begin{eqnarray}\label{L3.14}
{\Expec{\sum_{z\in\mathbb{Z}^d\cap[-\frac{L}{2},\frac{L}{2})^d}({\rm osc}_{B_R(z)}(\zeta^p))^2}}
&\lesssim&
\Expec{\zeta^{2p}}^{1-\frac{1}{p}}
\Expec{\Big(\sum_{z\in\mathbb{Z}^d\cap[-\frac{L}{2},\frac{L}{2})^d}({\rm osc}_{B_R(z)}\zeta)^2\Big)^p}^\frac{1}{p}\nonumber \\
&&+\Expec{\Big(\sum_{z\in\mathbb{Z}^d\cap[-\frac{L}{2},\frac{L}{2})^d}({\rm osc}_{B_R(z)}\zeta)^2\Big)^p}.
\end{eqnarray}
The latter can be seen as follows: From the elementary real-variable estimate
\begin{equation}\nonumber
|\tilde\zeta^p-\zeta^p|\lesssim|\zeta|^{p-1}|\tilde\zeta-\zeta|+|\tilde\zeta-\zeta|^p,
\end{equation}
we obtain by definition of ${\rm osc}$ that
\begin{equation}\nonumber
{\rm osc}_{B_R(z)}(\zeta^p)\lesssim
|\zeta|^{p-1}{\rm osc}_{B_R(z)}\zeta+({\rm osc}_{B_R(z)}\zeta)^p.
\end{equation}
Using that the discrete $\ell^{2p}(\mathbb{Z}^d)$-norm is estimated by
the discrete $\ell^{2}(\mathbb{Z}^d)$-norm, this implies
\begin{eqnarray*}\nonumber
\lefteqn{\sum_{z\in\mathbb{Z}^d\cap[-\frac{L}{2},\frac{L}{2})^d}({\rm osc}_{B_R(z)}(\zeta^p))^2}\\
&\lesssim&
\zeta^{2(p-1)}\sum_{z\in\mathbb{Z}^d\cap[-\frac{L}{2},\frac{L}{2})^d}({\rm osc}_{B_R(z)}\zeta)^2
+\Big(\sum_{z\in\mathbb{Z}^d\cap[-\frac{L}{2},\frac{L}{2})^d}({\rm osc}_{B_R(z)}\zeta)^2\Big)^p.
\end{eqnarray*}
H\"older's inequality w.\ r.\ t.\ to $\Expec{\cdot}$ applied to the first r.\ h.\ s.\
term with exponents $(\frac{p}{p-1},p)$ yields (\ref{L3.14}).

\medskip

{\bf Step 2}. Conclusion in case of $p\ge 2$ (the other case is easier and not needed later).
It remains to treat the first r.\ h.\ s.\ term of (\ref{L3.11}).
By H\"older's inequality w.\ r.\ t.\  $\Expec{\cdot}$ we have
\begin{equation}\label{L3.30}
\Expec{\zeta^p}\le\Expec{\zeta^{2p}}^\frac{p-2}{2p-2}\Expec{\zeta^2}^\frac{p}{2p-2}.
\end{equation}
Using $\Expec{\zeta}=0$ we obtain from the original Spectral Gap Estimate applied to
$\zeta$ itself
\begin{equation}\label{L3.31}
\Expec{\zeta^2}\lesssim
\Expec{\sum_{z\in\mathbb{Z}^d\cap[-\frac{L}{2},\frac{L}{2})^d}({\rm osc}_{B_R(z)}\zeta)^2}
\stackrel{\mbox{\small Jensen}}{\le}
\Expec{\Big(\sum_{z\in\mathbb{Z}^d\cap[-\frac{L}{2},\frac{L}{2})^d}
({\rm osc}_{B_R(z)}\zeta)^2\Big)^p}^\frac{1}{p}.
\end{equation}
Inserting (\ref{L3.31}) into (\ref{L3.30}) we obtain
\begin{equation}\nonumber
\Expec{\zeta^p}^2
\lesssim\Expec{\zeta^{2p}}^\frac{p-2}{p-1}
\Expec{\Big(\sum_{z\in\mathbb{Z}^d\cap[-\frac{L}{2},\frac{L}{2})^d}({\rm osc}_{B_R(z)}\zeta)^2\Big)^p}^\frac{1}{p-1}.
\end{equation}
Inserting this into (\ref{L3.11}) and using Young's inequality yields the claim of the lemma.

\bigskip


{\sc Proof of Lemma \ref{L2}}.

\medskip

{\bf Step 1}. Regularity theory for $a$-harmonic functions.
We will use the following two ingredients from De Giorgi's theory
for uniformly elliptic equations: For any $a\in\Omega$ and any $a$-harmonic function $u$ in $B_2$ we have
\begin{equation}\label{L2.8}
\sup_{B_1}|u|\lesssim\Big(\int_{B_{2}}u^2\Big)^\frac{1}{2}.
\end{equation}
Moreover, there exists a H\"older exponent $\alpha>0$ only depending on $d$ and $\lambda$ such that
\begin{equation}\label{L2.14}
\sup_{x_1,x_2\in B_1}\frac{|u(x_1)-u(x_2)|}{|x_1-x_2|^\alpha}\lesssim\Big(\int_{B_{2}}|\nabla u|^2\Big)^\frac{1}{2}.
\end{equation}
Both ingredients follow from De Giorgi's theorem (see for instance \cite[Theorem~4.11]{Han-Lin-97}):
\begin{equation}\label{DeGiorgi}
\sup_{B_1} |u|+\sup_{x_1,x_2\in B_1}\frac{|u(x_1)-u(x_2)|}{|x_1-x_2|^\alpha}\lesssim\Big(\int_{B_{2}}u^2\Big)^\frac{1}{2}.
\end{equation}
Indeed, for \eqref{L2.14}, one may assume w.~l.~o.~g. that $\fint_{B_2}u=0$ so that
\eqref{L2.14} follows from \eqref{DeGiorgi} and Poincar\'e's inequality on $B_2$ for functions with mean value zero.
Here and in the sequel, we write $B_R=B_R(0)$ for brevity.
The crucial element of these estimates is that the constants depend on the coefficient field
$a$ only through the ellipticity ratio $\lambda$ (as indicated by the use of $\lesssim$).

\medskip

{\bf Step 2}. In this step, we derive an auxiliary a priori estimate involving dyadic annuli.
Let $u$ be a function and $g$ a vector field on the torus 
related by the elliptic equation $-\nabla\cdot a\nabla u=\nabla\cdot g$ and normalized by $\int_{[-\frac{L}{2},\frac{L}{2})^d}u=0$. 
We claim that if $g$ vanishes in $B_1$ we have
\begin{equation}\label{L2.2}
|u(0)|\lesssim\sum_{n=1}^\infty(2^n)^{1-\frac{d}{2}}\Big(\int_{B_{2^n}\setminus B_{2^{n-1}}}|g|^2\Big)^\frac{1}{2}.
\end{equation}
We note that this sum is actually finite since for $2^n\gg L$, the ball $B_{2^{n-1}}$
invades the entire torus so that the ``annulus'' $B_{2^n}\setminus B_{2^{n-1}}$ is actually void.
Estimate (\ref{L2.2}) will be derived from (\ref{L2.8}) and an elementary scaling argument.
Indeed, for $n\in\mathbb{N}$, we introduce
\begin{equation}\nonumber
g_n:=\left\{\begin{array}{ccc}g&\mbox{on}&B_{2^n}\setminus B_{2^{n-1}}\\
0&\mbox{else}&\end{array}\right\},
\end{equation}
so that $g=\sum_{n=1}^\infty g_n$. Let $u_n$ denote the solution of $-\nabla\cdot a\nabla u_n=
\nabla\cdot g_n$ on the torus normalized by $\int_{[-\frac{L}{2},\frac{L}{2})^d}u_n=0$,
so that $u=\sum_{n=1}^\infty u_n$. Hence by the triangle inequality,
for (\ref{L2.2}) it is enough to show
\begin{equation}\label{L2.21}
|u_n(0)|\lesssim(2^n)^{1-\frac{d}{2}}\Big(\int_{[-\frac{L}{2},\frac{L}{2})^d}|g_n|^2\Big)^\frac{1}{2}.
\end{equation}
We now give the argument for (\ref{L2.21}). Testing $-\nabla\cdot a\nabla u_n=\nabla\cdot g_n$
with $u_n$, using the uniform ellipticity, and Cauchy-Schwarz' inequality, we obtain
\begin{equation}\nonumber
\Big(\int_{[-\frac{L}{2},\frac{L}{2})^d}|\nabla u_n|^2\Big)^\frac{1}{2}\lesssim\Big(\int_{[-\frac{L}{2},\frac{L}{2})^d}|g_n|^2\Big)^\frac{1}{2}.
\end{equation}
Since $d>2$, Sobolev's embedding together with $\int_{[-\frac{L}{2},\frac{L}{2})^d}u_n=0$ yields
\begin{equation}\nonumber
\Big(\int_{[-\frac{L}{2},\frac{L}{2})^d}|u_n|^\frac{2d}{d-2}\Big)^\frac{d-2}{2d}\lesssim\Big(\int_{[-\frac{L}{2},\frac{L}{2})^d}|g_n|^2\Big)^\frac{1}{2}.
\end{equation}
By H\"older's inequality on $B_{2^{n-1}}$ with exponents $(\frac{d}{d-2},\frac{d}{2})$  we obtain
\begin{equation}\nonumber
\Big(\int_{B_{2^{n-1}}}|u_n|^{2}\Big)^\frac{1}{2}\lesssim 2^n\Big(\int_{[-\frac{L}{2},\frac{L}{2})^d}|g_n|^2\Big)^\frac{1}{2}.
\end{equation}
We note that $u_n$ is $a$-harmonic on $B_{2^{n-1}}$ (since $g_n$ vanishes there).
Hence by (\ref{L2.8}), which we rescale from
$B_2$ to $B_{2^{n-1}}$, we have
\begin{equation}\nonumber
|u_n(0)|\le \sup_{B_{2^{n-2}}}|u_n|\lesssim\Big((2^n)^{-d}\int_{B_{2^{n-1}}}|u_n|^2\Big)^\frac{1}{2}.
\end{equation}
The combination of the two last estimates yields (\ref{L2.21}).

\medskip

{\bf Step 3}. As a preliminary, we study the local dependence of $\nabla\phi+\xi$
on $a$: Let the two coefficient fields $a$ and $\tilde a$ agree outside $B_R$.
Then we have
\begin{equation}\label{L2.4}
\int_{B_R}|\nabla\phi(\tilde a;\cdot)+\xi|^2\lesssim\int_{B_R}|\nabla\phi(a;\cdot)+\xi|^2.
\end{equation}
Indeed, we note that the function $\phi(\tilde a;\cdot)-\phi(a;\cdot)$
satisfies
\begin{equation}\nonumber
-\nabla\cdot \tilde a\nabla(\phi(\tilde a;\cdot)-\phi(a;\cdot))=\nabla\cdot(\tilde a-a)(\nabla\phi(a;\cdot)+\xi).
\end{equation}
We test this equation with $\phi(\tilde a;\cdot)-\phi(a;\cdot)$ and obtain from uniform
ellipticity and Cauchy-Schwarz' inequality
\begin{equation}\nonumber
\int_{[-\frac{L}{2},\frac{L}{2})^d}|\nabla(\phi(\tilde a;\cdot)-\phi(a;\cdot))|^2
\lesssim \int_{[-\frac{L}{2},\frac{L}{2})^d}|(\tilde a-a)(\nabla\phi(a;\cdot)+\xi)|^2.
\end{equation}
Since by assumption, $\tilde a-a$ vanishes outside $B_R$, the above yields
\begin{eqnarray}\label{L2.25}
\lefteqn{\int_{B_R}|\nabla(\phi(\tilde a;\cdot)-\phi(a;\cdot))|^2}\nonumber\\
&\le&\int_{[-\frac{L}{2},\frac{L}{2})^d}|\nabla(\phi(\tilde a;\cdot)-\phi(a;\cdot))|^2
\lesssim \int_{B_R}|\nabla\phi(a;\cdot)+\xi|^2.
\end{eqnarray}
This implies (\ref{L2.4}) by the triangle inequality in $L^2(B_R)$.

\medskip

{\bf Step 4}. In this step, we derive the central deterministic estimate
\begin{eqnarray}\label{L2.13}
\lefteqn{\Big(\sum_{z\in\mathbb{Z}^d\cap[-\frac{L}{2},\frac{L}{2})^d\setminus B_{R+1}}({\rm osc}_{B_R(z)}
\phi(\cdot;0))^{2}\Big)^\frac{1}{2}}\nonumber\\
&\lesssim&\sum_{n=1}^\infty
(2^n)^{1-\frac{d}{2}}
\Big(\sum_{z\in\mathbb{Z}^d\cap B_{2^n+R}}\big(\int_{B_R(z)}|\nabla\phi+\xi|^2\big)^p\Big)^\frac{1}{2p}.
\end{eqnarray}
Given a coefficient field $a$ on the torus 
and a point on the integer lattice $z\in\mathbb{Z}^d\cap[-\frac{L}{2},\frac{L}{2})^d$, 
we denote by $a_z$ an arbitrary coefficient field on the torus
that agrees with $a$ outside $B_R(z)$. We note that the function $\phi(a_z;\cdot)-\phi(a;\cdot)$ 
satisfies
\begin{equation}\label{L2.3}
-\nabla\cdot a\nabla(\phi(a_z;\cdot)-\phi(a;\cdot))=\nabla\cdot(a_z-a)(\nabla\phi(a_z;\cdot)+\xi).
\end{equation}
Given a discrete field $\{\omega_z\}_{z\in\mathbb{Z}^d\cap[-\frac{L}{2},\frac{L}{2})^d}$ we consider the 
function $u$ and the vector field $g$ on the torus defined through
\begin{eqnarray*}\nonumber
u(x)&:=&\sum_{z\in\mathbb{Z}^d\cap[-\frac{L}{2},\frac{L}{2})^d}\omega_z(\phi(a_z;x)-\phi(a;x)),\\
g(x)&:=&\sum_{z\in\mathbb{Z}^d\cap[-\frac{L}{2},\frac{L}{2})^d}\omega_z(a_z(x)-a(x))(\phi(a_z;x)+\xi)
\end{eqnarray*}
and note that (\ref{L2.3}) translates into $-\nabla\cdot a\nabla u=\nabla\cdot g$.
Provided $\omega_z=0$ for $z\in B_{R+1}$, we have $g(x)=0$ for $x\in B_1$.
Under this assumption, we may apply (\ref{L2.2}) from Step 2 and obtain
\begin{eqnarray}
\lefteqn{\Big|\sum_{z\in\mathbb{Z}^d\cap[-\frac{L}{2},\frac{L}{2})^d}\omega_z(\phi(a_z;0)-\phi(a;0))\Big|}\nonumber\\
&\lesssim&\sum_{n=1}^\infty(2^n)^{1-\frac{d}{2}}\Big(\int_{B_{2^n}\setminus B_{2^{n-1}}}\Big|\sum_{z\in\mathbb{Z}^d\cap[-\frac{L}{2},\frac{L}{2})^d}\omega_z
(a_z-a)(\nabla\phi(a_z;\cdot)+\xi)\Big|^2\Big)^\frac{1}{2}.\nonumber
\end{eqnarray}
Since $|a_z-a|\le 1$ is supported in $B_R(z)$ and since
$\{B_R(z)\}_{z\in\mathbb{Z}^d}$ locally have a finite overlap, this turns into
\begin{eqnarray}
\lefteqn{\Big|\sum_{z\in\mathbb{Z}^d\cap[-\frac{L}{2},\frac{L}{2})^d}\omega_z(\phi(a_z;0)-\phi(a;0))\Big|}\nonumber\\
&\lesssim&\sum_{n=1}^\infty(2^n)^{1-\frac{d}{2}}\Big(\sum_{z\in\mathbb{Z}^d\cap[-\frac{L}{2},\frac{L}{2})^d}
\omega_z^2\int_{B_R(z)\cap B_{2^n}}|\nabla\phi(a_z,\cdot)+\xi|^2\Big)^\frac{1}{2}\nonumber\\
&\lesssim&\sum_{n=1}^\infty(2^n)^{1-\frac{d}{2}}\Big(\sum_{z\in\mathbb{Z}^d\cap B_{2^n+R}}
\omega_z^2\int_{B_R(z)}|\nabla\phi(a_z,\cdot)+\xi|^2\Big)^\frac{1}{2}\nonumber\\
&\stackrel{(\ref{L2.4})}{\lesssim}&\sum_{n=1}^\infty
(2^n)^{1-\frac{d}{2}}\Big(\sum_{z\in\mathbb{Z}^d\cap B_{2^n+R}}
\omega_z^2\int_{B_R(z)}|\nabla\phi(a,\cdot)+\xi|^2\Big)^\frac{1}{2}\nonumber\\
&\stackrel{\mbox{\small H\"older in $z$}}{\lesssim}&
\Big(\sum_{z\in\mathbb{Z}^d\cap[-\frac{L}{2},\frac{L}{2})^d}\omega_z^{2q}\Big)^\frac{1}{2q}\nonumber\\
&&\times\sum_{n=1}^\infty(2^n)^{1-\frac{d}{2}}
\Big(\sum_{z\in\mathbb{Z}^d\cap B_{2^n+R}}\big(\int_{B_R(z)}|\nabla\phi(a,\cdot)+\xi|^2\big)^p\Big)^\frac{1}{2p},\nonumber
\end{eqnarray}
where $p$ and $q$ are dual exponents, that is, $\frac{1}{p}+\frac{1}{q}=1$.
Since $\{\omega_z\}_{z\in\mathbb{Z}^d\cap[-\frac{L}{2},\frac{L}{2})^d}$ was arbitrary under the constraint that $\omega_z=0$ for
$z\in B_{R+1}$ this implies by the duality of $\ell^{2q}(\mathbb{Z}^d)$ and $\ell^\frac{2q}{2q-1}(\mathbb{Z}^d)$
\begin{eqnarray}
\lefteqn{\Big(\sum_{z\in\mathbb{Z}^d\cap[-\frac{L}{2},\frac{L}{2})^d\setminus B_{R+1}}
|\phi(a_z;0)-\phi(a;0)|^\frac{2q}{2q-1}\Big)^\frac{2q-1}{2q}}\nonumber\\
&\lesssim&\sum_{n=1}^\infty
(2^n)^{1-\frac{d}{2}}
\Big(\sum_{z\in\mathbb{Z}^d\cap B_{2^n+R}}\big(\int_{B_R(z)}|\nabla\phi(a,\cdot)+\xi|^2\big)^p\Big)^\frac{1}{2p}.\nonumber
\end{eqnarray}
Since for any $z\in\mathbb{Z}^d$, $a_z$ was an arbitrary coefficient field that agrees with $a$
outside $B_R(z)$, this implies by the definition of ${\rm osc}_{B_R(z)}$
\begin{eqnarray}
\lefteqn{\Big(\sum_{z\in\mathbb{Z}^d\cap[-\frac{L}{2},\frac{L}{2})^d\setminus B_{R+1}}({\rm osc}_{B_R(z)}
\phi(\cdot;0))^\frac{2q}{2q-1}\Big)^\frac{2q-1}{2q}}\nonumber\\
&\lesssim&\sum_{n=1}^\infty
(2^n)^{1-\frac{d}{2}}
\Big(\sum_{z\in\mathbb{Z}^d\cap B_{2^n+R}}\big(\int_{B_R(z)}|\nabla\phi+\xi|^2\big)^p\Big)^\frac{1}{2p}.\nonumber
\end{eqnarray}
On the l.\ h.\ s.\ we use that since $\frac{2q}{2q-1}\le 2$, the discrete $\ell^\frac{2q}{2q-1}(\mathbb{Z}^d)$-norm
dominates the discrete $\ell^2(\mathbb{Z}^d)$-norm to obtain (\ref{L2.13}).

\medskip

{\bf Step 5}. Using stationarity, we upgrade Step 4 to the stochastic estimate
\begin{equation}\label{L2.10}
\Expec{\Big(\sum_{z\in\mathbb{Z}^d\cap[-\frac{L}{2},\frac{L}{2})^d\setminus B_{R+1}}({\rm osc}_{B_R(z)}
\phi(\cdot;0))^{2}\Big)^p}^\frac{1}{2p}
\lesssim\Expec{\Big(\int_{B_R}|\nabla\phi+\xi|^2\Big)^p}^\frac{1}{2p}.
\end{equation}
Indeed, we start from (\ref{L2.13}) in Step 4
and apply the triangle inequality to the sum over $n$ w.\ r.\ t.\ the norm $L^{2p}(\Omega)$:
\begin{eqnarray}\label{L2.12}
\lefteqn{\Expec{\Big(\sum_{z\in\mathbb{Z}^d\cap[-\frac{L}{2},\frac{L}{2})^d\setminus B_{R+1}}{\rm osc}_{B_R(z)}
(\phi(\cdot;0))^{2}\Big)^p}^\frac{1}{2p}}\nonumber\\
&\lesssim&\sum_{n=1}^\infty
(2^n)^{1-\frac{d}{2}}
\bigg(\sum_{z\in\mathbb{Z}^d\cap B_{2^n+R}}
\Expec{\Big(\int_{B_R(z)}|\nabla\phi(a,\cdot)+\xi|^2\Big)^p}\bigg)^\frac{1}{2p}.
\end{eqnarray}
We now note that the stationarity (\ref{L2.11}) of $\phi$ also yields
\begin{equation}\nonumber
\nabla\phi(a;x+z)=
\nabla\phi(a(\cdot+z);x)
\end{equation}
and thus
\begin{equation}\nonumber
\int_{B_R(z)}|\nabla\phi(a,x')+\xi|^2 dx'
=\int_{B_R}|\nabla\phi(a(\cdot+z),x)+\xi|^2 dx.
\end{equation}
By stationarity of $\Expec{\cdot}$, cf.\ (\ref{stat}) applied to $\zeta(a)=\int_{B_R(z)}|\nabla\phi(a;x)+\xi|^2dx$, this implies
\begin{equation}\label{L2.31}
\Expec{\Big(\int_{B_R(z)}|\nabla\phi(\cdot,x')+\xi|^2 dx'\Big)^p}
=\Expec{\Big(\int_{B_R}|\nabla\phi(\cdot,x)+\xi|^2 dx\Big)^p}.
\end{equation}
Inserting this into (\ref{L2.12}) yields (because of $\sum_{z\in\mathbb{Z}^d\cap B_{2^n+R}}1\lesssim (2^n+R)^d$)
\begin{eqnarray}\label{L2.9}
\lefteqn{\Expec{\Big(\sum_{z\in\mathbb{Z}^d\cap[-\frac{L}{2},\frac{L}{2})^d\setminus B_{R+1}}{\rm osc}_{B_R(z)}
(\phi(\cdot;0))^{2}\Big)^p}^\frac{1}{2p}}\nonumber\\
&\lesssim&\sum_{n=1}^\infty
(2^n)^{1-\frac{d}{2}}(2^n+R)^\frac{d}{2p}\Expec{\Big(\int_{B_R}|\nabla\phi(a,\cdot)+\xi|^2\Big)^p}^\frac{1}{2p}.
\end{eqnarray}
Since for $p>\frac{d}{d-2}$ the exponent $1-\frac{d}{2}+\frac{d}{2p}<0$ is negative we have
$$\sum_{n=1}^\infty(2^n)^{1-\frac{d}{2}}(2^n+R)^\frac{d}{2p}\lesssim 1.$$
Hence (\ref{L2.9})
turns into the desired (\ref{L2.10}).

\medskip

{\bf Step 6}. It remains to treat $z\in\mathbb{Z}^d\cap B_{R+1}$ in (\ref{L2.b}). By stationarity,
it will be enough to consider $z=0$, cf.\ Step 7. In this step, we will derive
from the H\"older continuity a priori estimate (\ref{L2.14}) the deterministic estimate
\begin{equation}\label{L2.20}
{\rm osc}_{B_R}\phi(a;0)\lesssim \Big(\int_{B_{2R}}|\nabla\phi(a;\cdot)+\xi|^2\Big)^\frac{1}{2}.
\end{equation}
Let $a\in\Omega$ be given and $\tilde a\in\Omega$ agree with $a$ outside $B_R$ and
otherwise be arbitrary.
On the one hand, since $d>2$ and $\int_{[-\frac{L}{2},\frac{L}{2})^d}(\phi(\tilde a;\cdot)-\phi(a;\cdot))
\stackrel{(\ref{T.3})}{=}0$, we have by Sobolev's embedding
\begin{eqnarray}\label{L2.19}
\Big(\int_{B_R}|\phi(\tilde a;\cdot)-\phi(a;\cdot)|^\frac{2d}{d-2}\Big)^\frac{d-2}{2d}
&\le&\Big(\int_{[-\frac{L}{2},\frac{L}{2})^d}|\phi(\tilde a;\cdot)-\phi(a;\cdot)|^\frac{2d}{d-2}\Big)^\frac{d-2}{2d}\nonumber\\
&\lesssim&\Big(\int_{[-\frac{L}{2},\frac{L}{2})^d}|\nabla(\phi(\tilde a;\cdot)-\phi(a;\cdot))|^{2}\Big)^\frac{1}{2}\nonumber\\
&\stackrel{(\ref{L2.25})}{\lesssim}&\Big(\int_{B_R}|\nabla\phi(a;\cdot)+\xi|^2\Big)^\frac{1}{2}.
\end{eqnarray}
On the other hand, we obtain from (\ref{L2.14})
applied to the $a$-harmonic function $u(x)=\phi(a;x)+\xi\cdot x$ (rescaled from $B_1$ to $B_R$):
\begin{equation}\label{L2.16}
\sup_{x_1,x_2\in B_R}\frac{|\phi(a;x_1)-\phi(a;x_2)+\xi\cdot(x_1-x_2)|}{|x_1-x_2|^\alpha}
\lesssim\Big(\int_{B_{2R}}|\nabla\phi(a;\cdot)+\xi|^2\Big)^\frac{1}{2}.
\end{equation}
Replacing $a$ by $\tilde a$ in the above and using (\ref{L2.4}) from Step 3 (with $B_R$ replaced
by $B_{2R}$) we likewise have
\begin{equation}\label{L2.17}
\sup_{x_1,x_2\in B_R}\frac{|\phi(\tilde a;x_1)-\phi(\tilde a;x_2)+\xi\cdot(x_1-x_2)|}{|x_1-x_2|^\alpha}
\lesssim\Big(\int_{B_{2R}}|\nabla\phi(a;\cdot)+\xi|^2\Big)^\frac{1}{2}.
\end{equation}
Combining (\ref{L2.16}) and (\ref{L2.17}), we obtain
\begin{equation}\label{L2.18}
\sup_{x_1,x_2\in B_R}\frac{|(\phi(\tilde a;x_1)-\phi(a;x_1))-(\phi(\tilde a;x_2)-\phi(a;x_2))|}{|x_1-x_2|^\alpha}
\lesssim\Big(\int_{B_{2R}}|\nabla\phi(a;\cdot)+\xi|^2\Big)^\frac{1}{2}.
\end{equation}
By the following elementary interpolation estimate, valid for an arbitrary function $u$,
\begin{equation}\nonumber
\sup_{B_R}|u|\,\lesssim\,\Big(\int_{B_R}|u|^\frac{2d}{d-2}\Big)^\frac{d-2}{2d}
+\sup_{x_1,x_2\in B_R}\frac{|u(x_1)-u(x_2)|}{|x_1-x_2|^\alpha},
\end{equation}
we see that (\ref{L2.19}) and (\ref{L2.18}) combine to
\begin{equation}\nonumber
|\phi(\tilde a;0)-\phi(a;0)|\lesssim \Big(\int_{B_{2R}}|\nabla\phi(a;\cdot)+\xi|^2\Big)^\frac{1}{2}.
\end{equation}
Since $\tilde a$ was arbitrary besides agreeing with $a$ outside $B_R$, we obtain
(\ref{L2.20}) by definition of ${\rm osc}$.

\medskip

{\bf Step 7}. We upgrade Step 6 to the stochastic estimate
\begin{equation}\label{L2.30}
\Expec{\Big(\sum_{z\in\mathbb{Z}^d\cap B_{R+1}}({\rm osc}_{B_R(z)}\phi(\cdot;0))^2\Big)^p}
\lesssim \Expec{\Big(\int_{B_{3R+1}}|\nabla\phi+\xi|^2\Big)^p}.
\end{equation}
Indeed, (\ref{L2.20}) from Step 6, with the origin replaced by $z$, implies after
summation
\begin{equation}\nonumber
\sum_{z\in\mathbb{Z}^d\cap B_{R+1}}({\rm osc}_{B_R(z)}\phi(\cdot;0))^2
\lesssim \int_{B_{3R+1}}|\nabla\phi+\xi|^2.
\end{equation}
Taking the $p$-th power and the expectation yields (\ref{L2.30}).

\medskip

{\bf Step 8}. From Steps 5 and 7 we learn that (\ref{L2.b}) is satisfied with $B_1$ replaced
by $B_{3R+1}$ on the r.\ h.\ s.\ . We appeal once more to stationarity to get for a generic $R\lesssim 1$
\begin{equation}\label{L2.32}
\Expec{\Big(\int_{B_R}|\nabla\phi+\xi|^2\Big)^p}\lesssim
\Expec{\Big(\int_{B_1}|\nabla\phi+\xi|^2\Big)^p}.
\end{equation}
Indeed, there exist points $z_1,\cdots,z_N$ on the torus such that
$B_R\subset \bigcup_{n=1}^NB_1(z_n)$ and we can arrange for $N\lesssim 1$ because of $R\lesssim 1$. Thus we have
\begin{equation}\nonumber
\int_{B_R}|\nabla\phi+\xi|^2\le\sum_{n=1}^N\int_{B_1(z_n)}|\nabla\phi+\xi|^2.
\end{equation}
Taking the $p$-th power gives
\begin{equation}\nonumber
\Big(\int_{B_R}|\nabla\phi+\xi|^2\Big)^p\lesssim\sum_{n=1}^N\Big(\int_{B_1(z_n)}|\nabla\phi+\xi|^2\Big)^p;
\end{equation}
taking the expectation yields
\begin{equation}\nonumber
\Expec{\Big(\int_{B_R}|\nabla\phi+\xi|^2\Big)^p}\le\max_{n=1,\cdots,N}\Expec{\Big(\int_{B_1(z_n)}|\nabla\phi+\xi|^2\Big)^p}.
\end{equation}
By stationarity, cf.\ (\ref{L2.31}), this yields (\ref{L2.32}).

\bigskip


{\sc Proof of Lemma \ref{L1}}

{\bf Step 1}. We start by establishing the deterministic estimate 
\begin{equation}\label{L1.1}
\Big(\int_{B_1}|\nabla\phi+\xi|^2\Big)^p\lesssim\int_{B_2}(\phi+\xi\cdot x)^{2(p-1)}|\nabla\phi+\xi|^2,
\end{equation}
which we will use in form of 
\begin{equation}\label{L1.4}
\Big(\int_{B_1}|\nabla\phi+\xi|^2\Big)^p\lesssim\int_{B_2}(\phi^{2(p-1)}+1)(|\nabla\phi|^2+1).
\end{equation}
Estimate (\ref{L1.1}) relies on the fact that $u(x):=\phi(x)+\xi\cdot x$ is $a$-harmonic, that is,
\begin{equation}\nonumber
-\nabla\cdot a\nabla u=0.
\end{equation}
We test this equation with $\eta^2u$, where $\eta$ is a cut-off function for $B_1$ in $B_2$.
By uniform ellipticity we obtain
\begin{equation}\nonumber
\lambda\int(\eta|\nabla u|)^2\le 2\int|\nabla\eta||u|\,\eta|\nabla u|.
\end{equation}
We now use Young's inequality (and $p\ge 2 >1$) on the r.\ h.\ s.\ integrand in form of
\begin{equation}\nonumber
\frac{2}{\lambda}|\nabla\eta||u|\,\eta|\nabla u|
\le \frac{1}{2}(\eta|\nabla u|)^2+C(|\nabla\eta|u)^{2\frac{p-1}{p}}(\eta|\nabla u|)^\frac{2}{p},
\end{equation}
which yields
\begin{equation}\nonumber
\int(\eta|\nabla u|)^2\lesssim \int(|\nabla\eta|u)^{2\frac{p-1}{p}}(\eta|\nabla u|)^\frac{2}{p}.
\end{equation}
By the choice of $\eta$, this implies
\begin{equation}\nonumber
\int_{B_1}|\nabla u|^2\lesssim \int_{B_2}u^{2\frac{p-1}{p}}|\nabla u|^\frac{2}{p}.
\end{equation}
It remains to apply Jensen's inequality on the r.\ h.\ s.\ to obtain as desired
\begin{equation}\nonumber
\int_{B_1}|\nabla u|^2\lesssim \Big(\int_{B_2}u^{2(p-1)}|\nabla u|^2\Big)^\frac{1}{p}.
\end{equation}

\medskip

{\bf Step 2}. We continue with the deterministic estimate
\begin{equation}\label{L1.3}
\int_{[-\frac{L}{2},\frac{L}{2})^d}\phi^{2(p-1)}|\nabla\phi|^2\lesssim\int_{[-\frac{L}{2},\frac{L}{2})^d}\phi^{2(p-1)},
\end{equation}
which we will use in form of 
\begin{equation}\label{L1.5}
\int_{[-\frac{L}{2},\frac{L}{2})^d}(\phi^{2(p-1)}+1)|\nabla\phi|^2\lesssim\int_{[-\frac{L}{2}\frac{L}{2})^d}(\phi^{2(p-1)}+1)
\end{equation}
that follows from the combination of \eqref{L1.3} once with the generic exponent $p$
and once with the exponent $p=2$. 
Indeed, we test $-\nabla\cdot a(\nabla\phi+\xi)=0$ with the monotone-in-$\phi$ expression $\frac{1}{2p-1}\phi|\phi|^{2(p-1)}$
over the entire torus. Because of $\nabla\frac{1}{2p-1}\phi|\phi|^{2(p-1)}=\phi^{2(p-1)}\nabla\phi$
and by uniform ellipticity, we obtain
\begin{equation}\nonumber
\lambda\int_{[-\frac{L}{2},\frac{L}{2})^d}\phi^{2(p-1)}|\nabla\phi|^2\le\int_{[-\frac{L}{2},\frac{L}{2})^d}\phi^{2(p-1)}|\nabla\phi|.
\end{equation}
Using Cauchy-Schwarz' inequality on the r.\ h.\ s.\ of that inequality yields (\ref{L1.3}).

\medskip

{\bf Step 3}. Conclusion using stationarity. We take the expectation of (\ref{L1.4}):
\begin{equation}\nonumber
\Expec{\Big(\int_{B_1}|\nabla\phi+\xi|^2\Big)^p}\lesssim\Expec{\int_{B_2}(\phi^{2(p-1)}+1)(|\nabla\phi|^2+1)}.
\end{equation}
By stationarity, we have
\begin{equation}\nonumber
\Expec{\int_{B_2}(\phi^{2(p-1)}+1)(|\nabla\phi|^2+1)}
=|B_2|L^{-d}\Expec{\int_{[-\frac{L}{2},\frac{L}{2})^d}(\phi^{2(p-1)}+1)(|\nabla\phi|^2+1)}.
\end{equation}
We now use the expectation of (\ref{L1.5}):
\begin{equation}\nonumber
|B_2|L^{-d}\Expec{\int_{[-\frac{L}{2},\frac{L}{2})^d}(\phi^{2(p-1)}+1)|\nabla\phi|^2}
\lesssim L^{-d}\Expec{\int_{[-\frac{L}{2},\frac{L}{2})^d}(\phi^{2(p-1)}+1)}.
\end{equation}
We use once more stationarity in form of
\begin{equation}\nonumber
L^{-d}\Expec{\int_{[-\frac{L}{2},\frac{L}{2})^d}(\phi^{2(p-1)}+1)}=\Expec{\phi^{2(p-1)}}+1.
\end{equation}

\bigskip


{\sc Proof of Proposition \ref{P}}.

By Jensen's inequality, it is enough to prove the statement for $p>\frac{d}{d-2}$.
We apply Lemma \ref{L3} to $\zeta(a)=\phi(a;0)$. We note that by stationarity of $\phi$ and $\Expec{\cdot}$
we have
\begin{equation}\nonumber
\Expec{\phi}=\Expec{ L^{-d}\int_{[-\frac{L}{2},\frac{L}{2})^d}\phi}\stackrel{(\ref{T.3})}{=}0.
\end{equation}
Hence the statement of Lemma \ref{L3} assumes the form 
\begin{equation}\nonumber
\Expec{\phi^{2p}}
\lesssim\Expec{\Big(\sum_{z\in\mathbb{Z}^d\cap[-\frac{L}{2},\frac{L}{2})^d}({\rm osc}_{B_R(z)}\phi(\cdot;0))^2\Big)^p}.
\end{equation}
Estimating the r.\ h.\ s.\ by Lemmas \ref{L2} and \ref{L1}, this turns into
\begin{equation}\nonumber
\Expec{\phi^{2p}}
\le C\Big(\Expec{\phi^{2(p-1)}}+1\Big).
\end{equation}
We conclude by using Jensen's and Young's inequalities in form of 
$C\Expec{\phi^{2(p-1)}}\le C\Expec{\phi^{2p}}^\frac{p-1}{p}\le\frac{1}{2}\Expec{\phi^{2p}}+\tilde C$.

\bigskip


{\sc Proof of Theorem \ref{T}}.

{\bf Step 1}. 
Application of Lemma \ref{L4} to $\zeta=\xi'\cdot a_{hom}^L\xi$ yields
\begin{equation}\label{T.6}
\Expec{(\xi'\cdot a_{hom}^L\xi-\Expec{\xi'\cdot a_{hom}^L\xi})^2}\lesssim
\Expec{\sum_{z\in\mathbb{Z}^d\cap[-\frac{L}{2},\frac{L}{2})^d}({\rm osc}_{B_R(z)}\xi'\cdot a_{hom}^L\xi)^2}.
\end{equation}

\medskip

{\bf Step 2}. Deterministic estimate of the oscillation.
We first rewrite $\xi'\cdot a_{hom}^L\xi$ as
\begin{equation}\label{Ant1}
\xi'\cdot a_{hom}^L\xi\,=\, L^{-d}\int_{[-\frac{L}{2},\frac{L}{2})^d} (\nabla \phi'+\xi')\cdot a (\nabla \phi+\xi),
\end{equation}
where $\phi'(a;x)$ is the corrector associated with the pointwise transpose field $a^t$ of $a$ and direction $\xi'$.
Indeed,  \eqref{Ant1} holds by \eqref{T.3} since $\phi'$ is periodic.
We claim
\begin{equation}\label{T.7}
{\rm osc}_{B_R(z)}\xi'\cdot a_{hom}^L\xi\lesssim L^{-d}\Big(\int_{B_R(z)}|\nabla\phi'+\xi'|^2\Big)^{\frac{1}{2}}\Big(\int_{B_R(z)}|\nabla\phi+\xi|^2\Big)^{\frac{1}{2}}.
\end{equation}
Indeed, consider two arbitrary coefficient fields $a_0,a_1\in\Omega$ that agree outside $B_R(z)$.
We write for abbreviation $\phi_i(x)=\phi(a_i;x)$, $\phi_i'(x)=\phi'(a_i,x)$, and $a_{hom,i}^L=a_{hom}^L(a_i)$ for $i=0,1$. By definition
of ${\rm osc}$ it is enough to show
\begin{equation}\label{T.1}
L^d|\xi'\cdot a_{hom,1}^L\xi-\xi'\cdot a_{hom,0}^L\xi|\lesssim\Big(\int_{B_R(z)}|\nabla\phi_0'+\xi'|^2\Big)^{\frac{1}{2}}\Big(\int_{B_R(z)}|\nabla\phi_0+\xi|^2\Big)^{\frac{1}{2}}.
\end{equation}
Indeed, we have by definition of $a_{hom}^L$ and of $\phi,\phi'$
\begin{eqnarray*}
\lefteqn{L^d(\xi'\cdot a_{hom,1}^L\xi-\xi'\cdot a_{hom,0}^L\xi)}\\
&\stackrel{(\ref{Ant1})}{=}&
\int_{[-\frac{L}{2},\frac{L}{2})^d}(\nabla\phi_1'+\xi')\cdot a_1(\nabla\phi_1+\xi)
-\int_{[-\frac{L}{2},\frac{L}{2})^d}(\nabla\phi_0'+\xi')\cdot a_0(\nabla\phi_0+\xi)\\
&=&
\int_{[-\frac{L}{2},\frac{L}{2})^d}\nabla(\phi_1'-\phi_0')\cdot a_1(\nabla\phi_1+\xi)
+\int_{[-\frac{L}{2},\frac{L}{2})^d}(\nabla\phi_0'+\xi')\cdot a_0\nabla(\phi_1-\phi_0)\\
&&+\int_{[-\frac{L}{2},\frac{L}{2})^d}(\nabla\phi_0'+\xi')\cdot (a_1-a_0)(\nabla\phi_1+\xi)\\
&=&
\int_{[-\frac{L}{2},\frac{L}{2})^d}\nabla(\phi_1'-\phi_0')\cdot a_1(\nabla\phi_1+\xi)
+\int_{[-\frac{L}{2},\frac{L}{2})^d}\nabla(\phi_1-\phi_0)\cdot a_0^t(\nabla\phi_0'+\xi')\\
&&+\int_{[-\frac{L}{2},\frac{L}{2})^d}(\nabla\phi_0'+\xi')\cdot (a_1-a_0)(\nabla\phi_1+\xi).
\end{eqnarray*}
Using the equation (\ref{T.3}) for $\phi_1$ and for $\phi_0'$, the first two r.~h.~s. terms vanish and this identity turns into
\begin{equation*}
{L^d(\xi'\cdot a_{hom,1}^L\xi-\xi'\cdot a_{hom,0}^L\xi)}\\
\,{=}\,\int_{[-\frac{L}{2},\frac{L}{2})^d}(\nabla\phi_0'+\xi')\cdot (a_1-a_0)(\nabla\phi_1+\xi),
\end{equation*}
so that we obtain
\begin{equation*}
{L^d|\xi'\cdot a_{hom,1}^L\xi-\xi'\cdot a_{hom,0}^L\xi|}
\,\le\,\Big(\int_{B_R(z)}|\nabla\phi_0'+\xi'|^2\int_{B_R(z)}|\nabla\phi_1+\xi|^2\Big)^\frac{1}{2}.
\end{equation*}
Now (\ref{T.1}) follows from this and Step 3 in the proof of Lemma \ref{L2}
in form of $\int_{B_R(z)}|\nabla\phi_1+\xi|^2\lesssim \int_{B_R(z)}|\nabla\phi_0+\xi|^2$.

\medskip

{\bf Step 3}. Stochastic estimate based on Proposition \ref{P}. We claim
\begin{equation}\label{T.5}
\Expec{\Big(\int_{B_R(z)}|\nabla\phi+\xi|^2\Big)^2}\lesssim 1.
\end{equation}
Indeed, by Step 8 from the proof of Lemma \ref{L2} (and stationarity to replace $z$ by $0$) we have 
\begin{equation}\nonumber
\Expec{\Big(\int_{B_R(z)}|\nabla\phi+\xi|^2\Big)^2}\lesssim 
\Expec{\Big(\int_{B_1}|\nabla\phi+\xi|^2\Big)^2}.
\end{equation}
An application of Lemma \ref{L1} with $p=2$ and of Proposition \ref{P} with $p=1$ yields (\ref{T.5}).
Since the ensemble $\expec{\cdot}'$ that is obtained from $\expec{\cdot}$ as pushforward
under $a\mapsto a^t$ satisfies our assumptions, we have
\begin{equation}
\Expec{\Big(\int_{B_R(z)}|\nabla\phi'+\xi'|^2\Big)^2} \,=\,\Expec{\Big(\int_{B_R(z)}|\nabla\phi+\xi'|^2\Big)^2}'\,\lesssim\, 1\label{T.5b}.
\end{equation} 

\medskip

{\bf Step 4}. Conclusion:
\begin{eqnarray*}
\lefteqn{\Expec{(\xi'\cdot a_{hom}^L\xi-\Expec{\xi'\cdot a_{hom}^L\xi})^2}}\\
&\stackrel{(\ref{T.6})}{\lesssim}&
\Expec{\sum_{z\in\mathbb{Z}^d\cap[-\frac{L}{2},\frac{L}{2})^d}({\rm osc}_{B_R(z)}\xi'\cdot a_{hom}^L\xi)^2}\\
&\stackrel{(\ref{T.7})}{\lesssim}&
L^{-2d}\Expec{\sum_{z\in\mathbb{Z}^d\cap[-\frac{L}{2},\frac{L}{2})^d}
\int_{B_R(z)}|\nabla\phi'+\xi'|^2\int_{B_R(z)}|\nabla\phi+\xi|^2}\\
&\leq & L^{-2d}\sum_{z\in\mathbb{Z}^d\cap[-\frac{L}{2},\frac{L}{2})^d}
\Expec{\Big(\int_{B_R(z)}|\nabla\phi'+\xi'|^2\Big)^2}^{\frac{1}{2}}\Expec{\Big(\int_{B_R(z)}|\nabla\phi+\xi|^2\Big)^2}^{\frac{1}{2}}\\
&\stackrel{(\ref{T.5}),(\ref{T.5b})}{\lesssim}&L^{-d}.
\end{eqnarray*}

\bigskip



\begin{thebibliography}{10}

\bibitem{Biskup-Salvi-Wolff-14}
M.~Biskup, M.~Salvi, and T.~Wolff.
\newblock {A central limit theorem for the effective conductance: I. Linear
  boundary data and small ellipticity constrasts}.
\newblock {\em Commun. Math. Phys.}, 2014.

\bibitem{Bourgeat-99}
A.~Bourgeat and A.~Piatnitski.
\newblock Estimates in probability of the residual between the random and the
  homogenized solutions of one-dimensional second-order operator.
\newblock {\em Asympt. Anal.}, 21(3-4):303--315, 1999.

\bibitem{Gloria-Neukamm-Otto-14}
A.~Gloria, S.~Neukamm, and F.~Otto.
\newblock Quantification of ergodicity in stochastic homogenization: optimal
  bounds via spectral gap on {G}lauber dynamics.
\newblock {\em Invent. Math.}, 2014.
\newblock DOI 10.1007/s00222-014-0518-z.

\bibitem{Gloria-Nolen-14}
A.~Gloria and J.~Nolen.
\newblock Quantitative central limit theorem for the effective conductance.
\newblock In preparation.

\bibitem{Gloria-Otto-10b}
A.~Gloria and F.~Otto.
\newblock Quantitative results on the corrector equation in stochastic
  homogenization of linear elliptic equations.
\newblock arXiv:1409.0801.

\bibitem{Gloria-Otto-09}
A.~Gloria and F.~Otto.
\newblock An optimal variance estimate in stochastic homogenization of discrete
  elliptic equations.
\newblock {\em Ann. Probab.}, 39(3):779--856, 2011.

\bibitem{Gloria-Otto-09b}
A.~Gloria and F.~Otto.
\newblock An optimal error estimate in stochastic homogenization of discrete
  elliptic equations.
\newblock {\em Ann. Appl. Probab.}, 22(1):1--28, 2012.

\bibitem{Han-Lin-97}
Q.~Han and F.~Lin.
\newblock {\em Elliptic partial differential equations}.
\newblock Courant Institute of Mathematical Sciences, New York, 1997.

\bibitem{Kozlov-79}
S.M. Kozlov.
\newblock The averaging of random operators.
\newblock {\em Mat. Sb. (N.S.)}, 109(151)(2):188--202, 327, 1979.

\bibitem{Marahrens-Otto-13}
D.~Marahrens and F.~Otto.
\newblock Annealed estimates on the {Green's} function.
\newblock MPI Preprint 69/2012.

\bibitem{Mourrat-Otto-14}
J.-C. Mourrat and F.~Otto.
\newblock Correlation structure of the corrector in stochastic homogenization.
\newblock  arXiv:1402.1924.

\bibitem{Naddaf-Spencer-98}
A.~Naddaf and T.~Spencer.
\newblock Estimates on the variance of some homogenization problems.
\newblock Preprint, 1998.

\bibitem{Nolen-11}
J.~Nolen.
\newblock Normal approximation for a random elliptic equation.
\newblock {\em Probab. Th. Rel. Fields}, 2013.

\bibitem{Papanicolaou-Varadhan-79}
G.C. Papanicolaou and S.R.S. Varadhan.
\newblock Boundary value problems with rapidly oscillating random coefficients.
\newblock In {\em Random fields, {V}ol. {I}, {II} ({E}sztergom, 1979)},
  volume~27 of {\em Colloq. Math. Soc. J\'anos Bolyai}, pages 835--873.
  North-Holland, Amsterdam, 1981.

\bibitem{Rossignol-12}
R.~Rossignol.
\newblock Noise-stability and central limit theorems for effective resistance
  of random electric networks.
\newblock arXiv:1206.3856.

\bibitem{Yurinskii-86}
V.~V. Yurinski{\u\i}.
\newblock Averaging of symmetric diffusion in random medium.
\newblock {\em Sibirskii Matematicheskii Zhurnal}, 27(4):167--180, 1986.

\end{thebibliography}
\end{document}